\documentclass[a4paper,10pt]{amsart}
\newif\ifpdf
\ifx\pdfoutput\undefined
    \pdffalse
\else
    \pdfoutput=1
    \pdftrue
\fi
\usepackage{enumerate}
\usepackage{amssymb}
\usepackage{mathrsfs}
\usepackage{amscd}
\usepackage[all]{xy}
\usepackage{graphicx}

\newcommand{\bc}{\mathbb C}

\newtheorem{thm0}{Theorem}
\newtheorem{cor0}{Corollary}

\newtheorem{prop0}{Proposition}

\theoremstyle{remark}

\theoremstyle{definition}

\theoremstyle{remark}

\newcommand\lineto{\ar@{-}}
\newcommand\dashto{\ar@{--}}
\newcommand\dotto{\ar@{.}}

\newcommand{\comment}[1]{}
\newbox\mybox
\def\overtag#1#2#3{\setbox\mybox\hbox{$#1$}\hbox to
  0pt{\vbox to 0pt{\vglue-#3\vglue-\ht\mybox\hbox to \wd\mybox
      {\hss$\ss#2$\hss}\vss}\hss}\box\mybox}
\def\undertag#1#2#3{\setbox\mybox\hbox{$#1$}\hbox to 0pt{\vbox to
    0pt{\vglue#3\vglue\ht\mybox\hbox to \wd\mybox
      {\hss$\ss#2$\hss}\vss}\hss}\box\mybox}
\def\lefttag#1#2#3{\hbox to 0pt{\vbox to 0pt{\vss\hbox to
      0pt{\hss$\ss#2$\hskip#3}\vss}}#1}
\def\righttag#1#2#3{\hbox to 0pt{\vbox to 0pt{\vss\hbox to
      0pt{\hskip#3$\ss#2$\hss}\vss}}#1}
\let\ss\scriptstyle

\def\Dot{\lower.2pc\hbox to 2.5pt{\hss$\bullet$\hss}}
\def\Circ{\lower.2pc\hbox to 2.5pt{\hss$\circ$\hss}}
\def\Vdots{\raise5pt\hbox{$\vdots$}}
\def\splicediag#1#2{\xymatrix@R=#1pt@C=#2pt@M=0pt@W=0pt@H=0pt}

\renewcommand\frame[2][3pt]{\hbox{$\vcenter{\hbox{\vrule\vbox
{\hrule\kern#1\hbox{\kern#1$#2$\kern#1}\kern#1\hrule}\vrule}}$}}

\newcommand{\ssbull}{\raise.2ex\hbox{${\scriptscriptstyle\bullet}$}}

\newcommand{\moplus}{\hbox{$\bigoplus$}}

\newcommand{\bC}{{\mathbb C}}

\newcommand{\bQ}{{\mathbb Q}}

\newcommand{\bZ}{{\mathbb Z}}

\newcommand{\cD}{{\mathcal D}}

\newcommand{\cO}{{\mathcal O}}

\newcommand{\Gr}{\text{{\rm Gr}}}
\newcommand{\Ker}{\text{{\rm Ker}}}

\def\mopls{\hbox{$\bigoplus$}}

\def\mcap{\hbox{$\bigcap$}}

\def\Ker{\operatorname{Ker}}

\def\ztop{Z^{}_{\operatorname{top}}(f,s)}
\def\ztopl{Z^{}_{\operatorname{top,0}}(f,s)}
\def\f{f^{-1}\{0\}}

\begin{document}

\title[Maximal poles of zeta functions, roots $b$-functions and  monodromy Jordan blocks]{On `maximal' poles of zeta functions, roots of $b$-functions and  monodromy Jordan blocks}
\author{A. Melle-Hern\'andez, T. Torrelli and Willem Veys}

\address{Facultad de Matem\'aticas\\ Universidad Complutense\\
Plaza de Ciencias~3\\ E-28040, Madrid, Spain}
\address{Laboratoire Jean Alexandre Dieudonn\'e\\
Universit\'e de Nice - Sophia Antipolis\\
Facult\'e des Sciences, Parc Valrose\\
06108 Nice Cedex 02, France}
\address{University of Leuven\\ Department of Mathematics\\ Celestijnenlaan 200 B,
B-3001 \\ Leuven (Heverlee), Belgium}

\email{amelle@mat.ucm.es, tristan.torrelli@laposte.net,
wim.veys@wis.kuleuven.be}

\thanks{\noindent The first author is partially supported by Spanish Contract MTM2007-67908-C02-02.
The third author is partially supported by the Fund of Scientific
Research - Flanders (G.0318.06).}



\dedicatory{Dedicated with admiration to C.T.C.~Wall on the
occasion of his seventieth birthday}

\begin{abstract}
The main objects of study in this paper are the poles of several
local zeta functions: the Igusa, topological and motivic zeta
function associated to a polynomial or (germ of) holomorphic
function in $n$ variables. We are interested in poles of maximal
possible order $n$. In all known cases (curves, non-degenerate
polynomials) there is at most one pole of maximal order $n$ which
is then given by the log canonical threshold of the function at
the corresponding singular point.

For an isolated singular point we prove that if the log canonical
threshold yields a pole of order $n$ of the corresponding (local)
zeta function, then it induces a root of the Bernstein-Sato
polynomial of the given function of multiplicity $n$ (proving one
of the cases of the strongest form of a conjecture of
Igusa-Denef-Loeser). For an arbitrary singular point we show under
the same assumption that the monodromy eigenvalue induced by the
pole has \lq a Jordan block of size $n$ on the (perverse) complex
of nearby cycles\rq.

\end{abstract}

\maketitle

\section*{Introduction}

\medskip \noindent {\bf 0.1.} Let $f : X \rightarrow \bC$ be a
non-constant analytic function on an open part $X$ of $\bC^n$. The
\lq classical\rq\ complex zeta function associated to $f$ is an
integral
$$Z_\varphi(f;s) := \int_X | f(x) |^{2s} \varphi(x) dx\wedge d\bar x$$
for $s \in \bC$ with $\Re (s) > 0$, where $\varphi$ is a
$C^\infty$ function with compact support on $X$. (Here and further
$x = (x_1,\cdots,x_n)$ and $dx = dx_1 \wedge \cdots \wedge dx_n$.)
One verifies that $Z_\varphi(f;s)$ is holomorphic in $s$. Either
by resolution of singularities \cite{BG}, or using the
Bernstein-Sato polynomial $b_f(s)$ of $f$ \cite{be:72}, one can
show that it admits a meromorphic continuation to $\bC$. The
second method also yields that each pole of $Z_\varphi(f;s)$ is a
translate by a nonpositive integer of a root of $b_f(s)$. And
moreover, for a root $s_0$ of $b_f(s)$, the order of $s_0 -m$ as
pole of $Z_\varphi(f;s)$ is at most the multiplicity of $s_0$ as
root of $b_f(s)$ \cite{Igusa-book}. In particular a pole of
(maximal) order $n$ induces a root of multiplicity $n$.

\bigskip
\noindent {\bf 0.2.} Let now $f : X \rightarrow \bQ_p$ be a
non-constant ($\bQ_p$-)analytic function on a compact open $X
\subset \bQ^n_p$, where $\bQ_p$ denotes the field of $p$-adic
numbers. Let $| \cdot |_p$ and $| dx |$ denote the $p$-adic norm
and the Haar measure on $\bQ^n_p$, normalized in the standard way.
The $p$-adic integral
$$Z_p(f;s) := \int_X | f(x)|^s_p |dx|,$$
again defined for $s \in \bC$ with $\Re(s) > 0$, is called the
($p$-adic) Igusa zeta function of $f$. Using resolution of
singularities Igusa \cite{ig:74},\cite{ig:74b} showed that it is a rational function
of $p^{-s}$; hence it also admits a meromorphic continuation to
$\bC$. (Everything can be generalized to finite extensions of
$\bQ_p$.) There are various `algebro-geometric' zeta functions,
related to the $p$-adic Igusa zeta functions: the motivic, Hodge
and topological zeta functions. We recall the definition of the
local and global version of the topological zeta function.

Let $f:(\bc^{n},0)\to(\bc,0)$ be a germ of a non-zero holomorphic
function $f$ (resp. $f\in \bc[x_1,\ldots,x_n]$ non-zero,
$f(0)=0$). Let $B$ be an open ball centered at the origin. Let
$\pi : X \to B$ (resp. $\pi : X \to \bc^{n}$) be an embedded
resolution of $(\f,0)$ (resp. $\f$).
We denote by $E_i, i\in J$, the irreducible components of
$\pi^{-1}(\f)_{\text{red}}$. Let $N_i$ (resp. $\nu_i-1$) be  the
multiplicity of $f \circ \pi$ (resp. of $\pi^*(dx_1\wedge\ldots
\wedge dx_n)$) at a generic point of $E_i$. For $I \subset J$, we
set $E_I := \cap_{i\in I} E_i$ and $E_I^{\circ}:= E_I \setminus
\left(\cup_{j\not\in I} E_j \right)$.


The {\it local topological zeta function} $\ztopl$ (resp. {\it
topological zeta function} $\ztop$) of $f$ at $0$ (resp. of $f$)
is the rational function defined by
$$
\ztopl:=\sum_{I\subset J}\chi(E_I^{\circ}\cap
\pi^{-1}\{0\})\prod_{i\in I}\frac{1}{\nu_i+N_i s}\in
\bQ(s),\leqno(*)
$$
$$
\ztop:=\sum_{I\subset J}\chi(E_I^{\circ})\prod_{i\in
I}\frac{1}{\nu_i+N_i s}\in \bQ(s),
$$
respectively. In \cite{dl:92}, Denef and Loeser proved that these
rational functions are well-defined (they do not depend on the
resolution $\pi$), by expressing them as a kind of limit of
$p$-adic Igusa zeta functions. We just mention that the motivic
zeta function specializes to the topological zeta function and to
the various $p$-adic Igusa zeta functions (for almost all $p$).

\medskip
In this paper we study a piece of a remarkable conjecture of
Igusa-Denef-Loeser, relating the poles of these zeta functions to
roots of the Bernstein-Sato polynomial, modeled on the result for
$Z_\varphi(f;s)$. We will treat poles of (maximal possible) order
$n$. For the topological zeta function it is clear that these
occur if and only if there exist $n$ different components $E_i$
with the same quotient $\nu_i/N_i$ and having a non-empty
intersection. For the other zeta functions the situation is
analogous. For that reason we formulate everything in terms of the
\lq simplest\rq\ zeta function, being the topological one. Our
results are however valid also for the other mentioned zeta
functions.

\medskip
{\bf Conjecture 1.} {\it The poles of $\ztopl$ are roots of the
local Bernstein-Sato polynomial $b_{f,0}(s).$}

\medskip
{\bf Conjecture 2.}  {\it The function $b_{f,0}(s)\cdot \ztopl$ is
a polynomial.}

\medskip
\noindent Conjecture 2 is a stronger version of Conjecture 1,
saying that the order of a pole $s_0$ of $\ztopl$ is at most the
multiplicity of $s_0$ as root of $b_{f,0}(s)$. For curves ($n=2$)
Conjecture 1 was proved by Loeser \cite{lo:88}. In that paper he
also verified Conjecture 2 for {\it reduced} $f$. For arbitrary
$n$ these conjectures are still wide open.

\bigskip
\noindent {\bf 0.3.} There is a well known relation between roots
of Bernstein-Sato polynomials and monodromy eigenvalues of $f$. In
particular, if $s_0$ is a root of $b_{f,0}(s)$, then $\exp(2\pi
is_0)$ is an eigenvalue of the monodromy acting on some cohomology
group of the (local) Milnor fibre of $f$ at some point of the germ
of $f^{-1}\{0\}$ at $0$ (equivalently; $\exp(2\pi is_0)$ is a
monodromy eigenvalue on the nearby cycle complex $\psi_f \bC$). So
the following conjecture, relating poles of $\ztopl$ to monodromy
eigenvalues, is implied by Conjecture 1.

\medskip
{\bf Conjecture 3.}  {\it If $s_0$ is a pole of $\ztopl$, then
$\exp(2\pi is_0)$ is an eigenvalue of the local monodromy acting
on some cohomology group of the Milnor fibre of $f$ at some point
of the germ of $f^{-1}\{0\}$ at $0$.}

\medskip

When $(f^{-1}\{0\},0)$ is a germ of an {\it isolated} singularity,
then a result of Varchenko \cite{varchenko} relates the
multiplicity of a root of $b_{f,0}(s)$ to the size of the
monodromy Jordan blocks for the associated monodromy eigenvalue. A
root of multiplicity $n$ corresponds essentially to a Jordan block
of size $n$ (see Theorem \ref{varchenko} for the precise
formulation).

This is certainly not true in general for non-isolated
singularities: for any homogeneous $f$ its monodromy is finite and
hence all Jordan blocks have size $1$. And for instance when
$f=\prod_{i=1}^n x_i^N$ we have that $b_{f,0}(s) = \prod_{j=1}^N
(s-j/N)^n$. The \lq right\rq\ generalization of Varchenko's result
should be stated in terms of the sub-complex $\psi_{f,\lambda}
\bC$ of the nearby cycle complex $ \psi_{f}\bC $; see Section 1.

\bigskip
\noindent {\bf 0.4.} In this paper we investigate for an arbitrary
$f$ in $n$ variables, assuming that its topological zeta function
has a pole $s_0$ of maximal order $n$, the implications concerning
$s_0$ being a root of $b_{f}(s)$ of multiplicity $n$, and
concerning a possible associated monodromy Jordan block of size
$n$. In a forthcoming paper we will study the case $n=2$ more in
detail, in particular for non-reduced $f.$

\bigskip
\noindent {\bf 0.5.} With the notation of 0.2 the  \emph{log
canonical threshold} $c_0(f)$ of $f$ at $0$ (resp. $c(f)$ of $f$)
is defined as
$$
c_0(f):=\min_{i\in J: 0\in \pi(E_i)}\{\nu_i/N_i\}, \qquad
c(f):=\min_{i\in J: \f\cap \pi(E_i)\ne \emptyset}\{\nu_i/N_i\},
$$
see e.g. Proposition 8.5 in \cite{kollar}. It does not depend on
the resolution $\pi$ since e.g. $-c(f)$ (resp. $-c_0(f)$) is the
root closest to the origin of the Bernstein-Sato polynomial
$b_{f}(s)$ (resp.  $b_{f,0}(s)$) of $f$ (at $0$), see e.g. Theorem 10.6
in \cite{kollar}. (In fact by results
of Lichtin and Kashiwara every root of $b_f(s)$ is of the form
$-\frac{\nu_i+k}{N_i}$, for some $i\in J$ and some integer $k\geq
0$, see Theorem 10.7 in \cite{kollar}.)

Clearly $-c_0(f)$ is the candidate pole of $\ztopl$ closest to the
origin. The third author formulated the following.

\medskip
{\bf Conjecture 4.}

{\it (1) $\ztopl$ has at most one pole of order $n$.

(2) If $\ztopl$ has in $s_0$ a pole of order $n$, then $s_0$ is
the pole closest to the origin of  $\ztopl$.}

\medskip
\noindent This conjecture is proved in case $n=2$ by himself
\cite{ve:95} and with Laeremans \cite{lv} when $f$ is
non-degenerate with respect to its Newton polyhedron and in
these cases $s_0 = -c_0(f)$ in (2).

Our main result  is roughly as
follows. Let $f:(\bc^{n},0)\to(\bc,0)$ be a germ of a holomorphic
function with $f(0)=0$ such that $s_0=-c_0(f)$ is a pole of order
$n$ of $\ztopl$. Denote $\lambda := \exp(2\pi i s_0)$. Then the
$\lambda$-characteristic subspace of the $(n-1)$-th cohomology of
the Milnor fibre of $f$ at $0$ has a nonzero $(2n-2$)-graded part
of its weight filtration. Morally \lq $\lambda$ has a Jordan block
of size $n$ on the perverse sheaf $\psi_{f}\bC$\rq. See
Theorem \ref{grn} and Corollary \ref{perverse} for a precise formulation.
The result of Varchenko implies then:

\begin{thm0}\label{bernstein}
Let $f:(\bc^{n},0)\to(\bc,0)$ be a germ of a non-zero holomorphic
function such that $(\f,0)$ is a germ of an isolated hypersurface
singularity. If $s_0=-c_0(f)$ is a pole of order $n$ of  $\ztopl$,
then $(s+c_0(f))^n$ divides the {Bernstein-Sato polynomial}
$b_{f,0}(s)$.

In such a case there exists an integer $N\geq 1$ such that $c_0(f)=1/N$ and either
\begin{itemize}
\item
$N=1$ and $(s+1)^n$ divides $b_{f,0}(s)$, or
\medskip
\item $N>1$ and
$(s+1/N)^n(s+2/N)^n\dots(s+(N-1)/N)^n(s+1)^n$ divides $b_{f,0}(s)$.
\end{itemize}
\end{thm0}





\medskip
\noindent If Conjecture 4 is true, then Theorems \ref{grn} and
\ref{bernstein} treat in fact \lq all\rq\ poles of maximal order
$n$. Our proof uses a result of Saito (Proposition \ref{key}) and
the ideas in the proof of the main result of van Doorn and
Steenbrink \cite{vDS}. In Theorem \ref{bernstein-g} we provide a
global version for polynomials $f$ and $\ztop$.

\medskip
\noindent {\it Acknowledgements.} The authors would like to thank
A.~Dimca, M.~Saito and J.~Steenbrink for useful discussions and
comments about the paper.

\medskip

\section{Preliminaries}
\medskip

\noindent {\bf 1.1. Monodromy.} Let $f$ be a  holomorphic function
on an $n$-dimensional complex manifold $X$.  Denote by $ X_t$ the
hypersurface $f^{-1}\{t\}$ for $t\in \bC$. Let $x\in X_0$ and
choose $\varepsilon,\eta>0$ with $\eta<<\varepsilon<<1$. The
restriction of $f$ to $\{ z\in X \,\big|\, |z-x|\leq\varepsilon,
\, 0<|f(z)|<\eta\}$ is a $C^\infty$ fibre bundle, {\it the Milnor
fibration}, whose typical fibre
$$
F_{f,x}:= \{z\in X \,\big|\, |z-x|\leq\varepsilon, f(z)=\delta\}
\quad\text{for}\,\,\, 0 < \delta <\eta
$$
is called the {\it Milnor fibre} of
$ f $ at $ x \in X_0$.
The Milnor fibre is endowed with the monodromy automorphism $M_{f,x}$
which induces an automorphism, denoted by $M_{f,x}^q$, on the cohomology groups
$ H^{q}(F_{f,x},\bC) $.

\medskip

Following Deligne \cite{deli} one has a sheaf theoretic version of
the previous constructions. Let $D$ be a small disk around the
origin in $\bC$, $D^*:=D\setminus\{0\}$ and ${\tilde D}^*$ the
universal covering of $D^*$. Consider the preimage $X^*$ of $D^*$
in $X$ and denote by ${\tilde X}^*$ the fibre product
$X^*\times_{D^*}{\tilde D}^*$. Let $i:X_0\to X$ be the inclusion
morphism and $j:{\tilde X}^*\to X$.


For the constructible sheaf $\bC_{X}$ on $X$ and for any $q\geq 0$, the
{\it nearby cycle sheaf} $R^q\psi_{f}\bc_{X}:=i^*R^qj_*j^*\bC_{X}$ is a constructible
sheaf on $X_0$. The deck transformation $(x,u)\mapsto (x,u+1)$ on
${\tilde  X}^*$ induces the action of a canonical monodromy automorphism
$T^q$ on $R^q\psi_{f}\bc_{X}$ such that the vector space
$(R^q\psi_{f}\bc_{X},T^q)_x$ with automorphism is canonically isomorphic to
$( H^{q}(F_{f,x},\bC), M_{f,x}^q) $.


In fact, working on the derived category of complexes with
automorphisms and bounded constructible cohomology, the {\it
nearby cycle complex} $\psi_{f}\bC_{X}$ on $X_0$ is defined by
$\psi_{f}\bC_{X}:=i^*R j_*j^*\bC_{X}$, see \cite{deli}. Recall
that the sheaf $\psi_{f}\bC_{X}[n-1]$ is a perverse sheaf. The
monodromy $T$ on the shifted perverse sheaf $\psi_{f}\bC_{X}$
 admits a decomposition $T=T_{s}T_u$ where $T_s$
is semisimple and $T_u$ is unipotent. For $ \lambda \in \bC $, let
$$
\psi_{f,\lambda}\bc_{X} = \Ker\,(T_{s}-\lambda) \subset
\psi_{f}\bC_{X}.\,\,\,
$$
\noindent There are also decompositions
$$
\psi_{f}\bC_{X} = \moplus_{\lambda}\psi_{f,\lambda}\bC_{X},\quad
H^{q}(F_{f,x},\bC) = \moplus_{\lambda}H^{q}(F_{f,x},\bC)_{\lambda},
$$
such that the action of $ T_{s} $ on $ \psi_{f,\lambda}\bC_{X}$
and on $H^{q}(F_{f,x},\bC)_{\lambda} $ is the multiplication by $
\lambda \in \bC^{*} $. The groups
$H^{q}(F_{f,x},\bC)_{\lambda}\oplus H^{q}(F_{f,x},\bC)_{\bar{\lambda}}={\mathcal
H}^q(\psi_{f,\lambda}\bC_{X})_x \oplus {\mathcal
H}^q(\psi_{f,\bar{\lambda}}\bC_{X})_x$ have a canonical Mixed Hodge
structure, see e.g.
\cite{[16]},\cite{[17]},\cite{st-kyoto}. Let $ W $ be
the weight filtration of this canonical mixed Hodge structure. The
following proposition is proved by M.~Saito (for a proof see
(1.1.3) and Proposition 1.7 in \cite{[36]}).

\medskip\noindent
\begin{prop0}\label{key} \cite{[36]}
Let $ N $ be the logarithm of the unipotent part $T_u$ of the
monodromy $ T $. If $ \Gr_{2n-2}^{W}H^{n-1}(F_{f,x},\bC)_{\lambda}
\ne 0 $, then $ N^{n-1} \ne 0 $ on $ \psi_{f,\lambda}\bC_{X} $ in
the category of shifted perverse sheaves.
\end{prop0}

Since $X$ is smooth and $n$-dimensional, $N^{n}=0$ on the nearby
cycle sheaf $ \psi_{f}\bC_{X}[n-1]$. This implies that the Jordan
blocks of the monodromy $M_{f,x}^q$ on the cohomology groups
$H^{q}(F_{f,x},\bC) $ have size $\leq q+1$ (see e.g. \cite{ds:03}
and references there). In fact, it is also proved there that the
support of the perverse sheaf $N^{n-1}  \psi_{f}\bC_{X}[n-1]$ is
empty or $0$-dimensional (see Proposition 0.5 in \cite{ds:03}).

\bigskip

\noindent {\bf 1.2. Bernstein-Sato polynomials.}
 Let $ X $ be a
complex $n$-dimensional manifold, resp. smooth algebraic variety,
and let $ X_0 $ be the hypersurface defined as the zero locus of a
holomorphic function, resp. regular function, $f$. Let $\cD_{X}$
be the ring of analytic, resp. algebraic, partial differential
operators associated to $X$.

The Bernstein-Sato polynomial (or $b$-function) $b_{f}(s)$ of $f$
is the unique monic polynomial of lowest degree satisfying
$$
b_{f}(s)f^{s} = Pf^{s+1}\quad\text{with}\,\,\,P\in \cD_{X}[s].
$$
It exists at least locally, and globally if $X$ is an affine
algebraic variety \cite{be:72},\cite{bj:79},\cite{sa:72}. Moreover
the $b$-function of a regular function $f$ and of its associated
analytic function coincide. Restricting to the stalk at a point $x
\in X_0$, one can also define the local $b$-function $b_{f,x}(s)$.
If $X$ is Stein, resp. affine, then $b_f(s)$ is the least common
multiple of these local $b$-functions.

Let $R_{f}$ be the set of the roots of $b_{f}(-s)$, and
$m_{\alpha}$ the multiplicity of $\alpha \in R_{f}$. Then $R_{f}
\subset \bQ_{>0}$, and $m_{\alpha} \le n$ because $b_{f}(s)$ is
closely related to the monodromy on the nearby cycle sheaf
$\psi_{f}\bC_{X}$, see e.g. \cite{[20]}. Set $\alpha_{f} =
\min R_{f}$; this number coincides with the log canonical
threshold, see \cite{kollar},\cite{[16]}.




\section{Monodromy on 
$\psi_{f}\bC_{X}$ and poles of zeta functions} 


\noindent {\bf 2.1.} We are interested in poles of maximal order
$n$ of $\ztopl$. Laeremans and the third author \cite{lv} proved
that every pole of maximal order of $\ztopl$ (or of $\ztop$) must
be of the form $-1/N$, for a positive integer  $N\geq 1$.

\medskip

\begin{thm0}\label{grn}
Let $f:(\bc^{n},0)\to(\bc,0)$ be a germ of a non-zero holomorphic
function. If $s_0=-c_0(f)$ is a pole of order $n$ of $\ztopl$,
then

$$\Gr_{2n-2}^{W}H^{n-1}(F_{f,0},\bC)_{\lambda} \ne 0 \quad  \text{ for } \lambda:= \exp(2\pi i s_0).$$

\noindent In such a case, there exists an integer $N\geq 1$, such that $c_0(f)=1/N$ and either
\begin{itemize}
\item
$N=1$ and $\Gr_{2n-2}^{W}H^{n-1}(F_{f,0},\bC)_{1} \ne 0 $ or
\medskip
\item $N>1$ and $\Gr_{2n-2}^{W}H^{n-1}(F_{f,0},\bC)_{\exp(2\pi i(-j/N))} \ne 0 $
for all $j$ with $1\leq j\leq N$.
\end{itemize}
\end{thm0}

\medskip

\begin{proof}
Assume $s_0$ is a pole of (maximal) order $n$ of $\ztopl$, then
write $s_0=-c_0(f)=-1/N$ for some integer $N\geq 1$ and  set
$\lambda:= \exp(2\pi i(-1/N))$. To show that
$\Gr_{2n-2}^{W}H^{n-1}(F_{f,0},\bC)_{\lambda} \ne 0 $ we will
adapt the proof of the main result of van Doorn and Steenbrink in
\cite{vDS}, see also Varchenko \cite{varchenko}.

\medskip

Let $B$ be an open ball centered at the origin. Let $\pi : X \to
B$ be an embedded resolution of the germ $(\f,0)$ which is an
isomorphism outside of the preimage of $\f$. Set $E:=\pi^{-1}(\f)$
and $E_x:=\pi^{-1}(0)$ and denote by $E_i, i\in J$, the
irreducible components of $E$. For $I \subset J$, put also $E_{I}
:= \mcap_{j\in I} E_{j}$.

\medskip

By the definition of the local topological zeta function, see
$(*)$, since $s_0$ is a pole of order $n$ of  $\ztopl$, there
exist $n$ irreducible components $E_0,\ldots,E_{n-1}$ of $E$ and
there is a point ${\tilde x}_0\in\cap_{i=0}^{n-1} E_i$ such that
$\pi({\tilde x}_0)=0$ and $\nu_i/N_i=1/N$ for all $0\leq i\le
n-1$.

We may assume that one of these irreducible components, called
e.g. $E_0$, is a K\"ahler compact non-singular variety. Otherwise
we blow up $X$ at ${\tilde x}_0$ and get a new configuration of
exceptional divisors where the new exceptional divisor is a
K\"ahler compact non-singular variety $E_0$ with $\nu_0/N_0=1/N$,
and we can choose a \lq new\rq\ ${\tilde x}_0$ on this $E_0$
satisfying the requirements above.

\medskip

To describe the quotient $ \Gr_{2n-2}^{W}H^{n-1}(F_{f,0},\bC)_{\lambda}$ one uses the fact that it is pure of type $(n-1,n-1)$ and therefore a quotient of the piece $F^{n-1}$ of the Hodge filtration.
These terms can be computed using the relative logarithmic de Rham complex.

Let $e$ be a common multiple of all multiplicities $N_i, i\in J$,
and let $\tilde \bc$ be another copy of $\bc$. Let $\tilde Y$ be
the normalization of the space $Y$ obtained from $X$ by the base
change $\sigma:{\tilde \bc}\to \bc: \sigma({\tilde t})={\tilde
t}^e$. Let $\rho:\tilde Y \to X$ and $\tilde f:\tilde Y\to {\tilde
\bc}$ be the natural projections maps. Let $D_i:=\rho^{-1}(E_i),
i\in J,$ and set $D:=\rho^{-1}(E)$, this is nothing but
$D=\cup_{i\in J} D_i$. Let $D_x:=\rho^{-1}(E_x)$. For every
$I\subset J$, let $D_{I}:=\rho^{-1}(E_I)$. The map $D_I\to E_I$ is
a cyclic cover of degree $\gcd(N_i,\, i\in I)$.

\[
\xymatrix{
\tilde Y \ar@/_/[ddr]_{\tilde f} \ar@/^/[drr]^\rho \ar@{.>}[dr]|{}\\
&Y={\tilde \bc} \times_\bc X\ar[d] \ar[r] & X\ar[d]_{f\circ \pi} \\
&{\tilde \bc} \ar[r]^{\sigma} &\bc}
\]

\medskip

\noindent By the semi-stable reduction theorem, $\pi$ and $e$ can be chosen in such a way that
$\tilde Y$ is smooth. The divisor $D={\tilde f}^{-1}(0)$ is a reduced normal crossing divisor, see \cite{st-kyoto}. From \cite{st-oslo}, see also \cite{st-kyoto}, there is an isomorphism

$$H^{q}(F_{f,0},\bC)\simeq  {\mathbb H}^q(D_x,\Omega^{\bullet}_{{\tilde Y}/{\tilde \bc}}(\log D)\otimes \cO_{D_x}),$$
so in particular $ \Gr_{F}^{p} H^{q}(F_{f,0},\bC)\simeq H^{q-p}(D_x,\Omega^{p}_{{\tilde Y}/{\tilde \bc}}(\log D)\otimes \cO_{D_x})$.
Then

\begin{multline*}
F^{n-1} H^{n-1}(F_{f,0},\bC)\simeq \Gr_{F}^{n-1} H^{n-1}(F_{f,0},\bC)\simeq   H^{0}(D_x,\Omega^{n-1}_{{\tilde Y}/{\tilde \bc}}
(\log\, D)\otimes \cO_{D_x})\\
\simeq   H^{0}(D_x,\Omega^{n}_{{\tilde Y}}(\log\, D)\otimes \cO_{D_x}).
\end{multline*}

\medskip

The following results can be deduced from Section 4 in \cite{varchenko}.
For every $\omega\in H^0(B,\Omega^n)\,$ define the {\it geometrical weight}
$g(\omega)$   with respect to the resolution $\pi$ as
$$
g(\omega):=\min_{i\in J} \left\{\frac{\text{ord}_{E_i}(\omega)+1}{N_i}\right\}.
$$

\noindent For every $\omega\in H^0(B,\Omega^n)\,$ with {\it
geometrical weight} $g(\omega)\leq 1$, define $R(\omega):={\tilde
f}^{\frac{-e}{N}} (\pi\rho)^*(\omega)$. Then $R(\omega)\in
H^{0}(\tilde Y,\Omega^{n}_{{\tilde Y}}(\log D))$. Let
$\sigma(\omega)$ be its Poincar\'e residue along $D_x$, that is
$\sigma(\omega)$ is the restriction to $D_x$ of ${\tilde
f}^{\frac{-e}{N}+1} (\pi\rho)^*(\omega)/d{\tilde f}$. Then
$\sigma(\omega)$ is an element in  ${F}^{n-1}
H^{n-1}(F_{f,0},\bC)$ and the semisimple part of the monodromy
acts on $\sigma(\omega)$ as
$$T_s(\sigma(\omega))=\exp(-2\pi ig(\omega)) \sigma(\omega).$$
The form $R(\omega)$ has a first order pole along $D_i$ if and only if $g(\omega)=\frac{\text{ord}_{E_i}(\omega)+1}{N_i}$, and else $R(\omega)$ is regular along $D_i$.

Consider the differential form $\eta=dx_1\wedge\ldots \wedge dx_n$.
Since $s_0=-c_0(f)=-1/N$ does not depend on the resolution,
$g(\eta)=1/N\leq 1$ does not depend on $\pi.$
In fact $1/N$ is the minimum and then $R(\eta)$ has a first order pole along $D_i$ if and only if $1/N=\nu_i/N_i$ and else $R(\eta)$ is regular along $D_i$.

\medskip

Let $D_{00}$ be one of the irreducible components of $D_0$. On the
open subspace $D_{00}^{\circ}= D_{00} \setminus \left(\cup_{j\ne
0} D_j \right)$, the restriction  of the map $\rho$ from
$D_{00}^{\circ}$ to $E_0^{\circ}$ is an \'etale cover. The form
$\sigma(\eta)$ is not equal to zero on $D_{00}^{\circ}$
and is in fact a meromorphic $(n-1)$-form with logarithmic poles
on $D_{00}\setminus D_{00}^{\circ}$.  Thus $R(\eta)$ defines a
non-zero element in $H^0({\tilde Y},\Omega^n_{\tilde
Y}(\log\,D))$. Notice that, according to Deligne's theorem
\cite{de:71},\cite{st-oslo},\cite{varchenko}, the class of
$\sigma(\eta)$ in $H^{n-1}(D_{00}^{\circ},\bc)$ is a
\emph{non-zero element} since $D_{00}$ is a projective manifold
and $D_{00}\setminus D_{00}^{\circ}\subset D_{00}$ is a divisor
with normal crossings.

\medskip

On an adequate chart on $\tilde Y$, with local coordinates
$y_1,\ldots,y_n$, the function ${\tilde f}$ is given by
$$\tilde f(y_1,\ldots, y_n)=\,y_1y_2\ldots y_{k}.$$
Moreover, if $y_1=0$ is the equation of the divisor $D_{00}$ (and
hence of $D_{0}$) and $y_j=0$ gives the divisor $D_{00}\cap D_j$, then
$$\sigma(\eta)=q(y_2,\ldots,y_n)\,\cdot y_2^{a_2}\dots y_k^{a_k}\frac{dy_2}{y_2}\wedge\ldots\wedge\frac{dy_k}{y_k}\wedge dy_{k+1}\wedge \ldots \wedge dy_n,$$
where $q(y_2,\ldots,y_n)$ is holomorphic, $q(0)\ne 0$ and $a_j=e\left(\frac{\nu_{j(i)}}{N_{j(i)}}-\frac{1}{N}\right)$.

In particular, at a preimage  $P_0$ of ${\tilde x}_0$ in $\tilde
Y$, $R(\eta)$ can be written locally as
$u\frac{dy_1}{y_1}\wedge\ldots\wedge\frac{dy_n}{y_n}$, with
$u(0)\ne 0$, because of the minimum. Considering for each $n$-fold
point $P$ on $D_x$ the multiple residue map Res$_P:H^0({\tilde
Y},\Omega^n_{\tilde Y}(\log\,D))\to \bc$, we have in particular
that Res$_{P_0}$ is surjective.

\medskip

Let $V_\lambda$ be the set of the $n$-fold points of $D$ that are
preimages of those $n$-fold points in $\cup_I E_I$ for which
$|I|=n$, $\lambda^{N_i}=1$ for all $i\in I$, and at least one
$E_i, i\in I,$ is an irreducible component of $E_x$. Then
$$ \Gr_{2n-2}^{W}H^{n-1}(F_{f,0},\bC)_{\lambda}\cong \text{Image} (\mopls_{P\in V_\lambda}
Res_{P}:H^0({\tilde Y},\Omega^n_{\tilde Y}(log\,D))\to
\bc^{V_\lambda});
$$
for the isomorphism see \cite{st-oslo} and \cite{vDS}.

Since $P_0\in V_\lambda$ and Res$_{P_0}$ is surjective then
$\Gr_{2n-2}^{W}H^{n-1}(F_{f,0},\bC)_{\lambda}\ne 0$ which concludes the proof.

\medskip

To show that $\Gr_{2n-2}^{W}H^{n-1}(F_{f,0},\bC)_{\exp(2\pi
i(-j/N))} \ne 0 $ also for $2\leq j\leq N$ we argue as follows,
see \cite{vDS}. Since the weight filtration is defined over $\bQ$,
it has a complex conjugation compatible with the one of
$\bc^{V_{\lambda}}$. Let $A_{P_0}=\{g(\omega):\, \omega\in
H^0(B,\Omega^n),\,\, g(\omega)\leq 1 \, \text{and}\,
\text{Res}_{P_0}(R(\omega))\ne 0\,\}$. Then
$g(\eta)=\frac{1}{N}\in A_{P_0}$. If $1/N<1$, then the complex
conjugate of $\text{Res}_{P_0}(R(\eta))$ in $\bc^{V_\lambda}$ is
an eigenvector of the semisimple part of the monodromy $T_s$ for
the eigenvalue  ${\bar \lambda}= \exp(2\pi i(1/N))$. In particular
there exists ${\tilde \eta}\in H^0 (B,\Omega^n)$ such that
$g({\tilde \eta})=\frac{N-1}{N}\in A_{P_0}$. After the remark on
page 230 in  \cite{vDS} one can prove that
$\{1/N,2/N,\ldots,(N-1)/N,1\}\subset A_{P_0}$. \end{proof}

\medskip

\noindent {\bf 2.2.} Using Proposition \ref{key} one `morally' obtains a Jordan block
of size $n$ in the category of shifted perverse sheaves.

\begin{cor0} \label{perverse} Let $f:(\bc^{n},0)\to(\bc,0)$ be a germ of a non-zero holomorphic function.
If $s_0=-c_0(f)$ is a pole of order $n$ of  $\ztopl$ and $\lambda:= \exp(2\pi i s_0)$, then $ N^{n-1} \ne 0 $ on $ \psi_{f,\lambda}\bC_{X} $ in
the category of shifted perverse sheaves.
\end{cor0}

\section{Applications for isolated hypersurface singularities} 

\noindent {\bf 3.1.} Let$f:(\bc^{n},0)\to(\bc,0)$ be a germ of a holomorphic function such
that $(\f,0)$ is a germ of an isolated hypersurface singularity.
The following result, by Varchenko \cite{varchenko}{ Theorem 1.4},
relates roots of the Bernstein-Sato polynomial ${b}_{f,0}(s)$ and
Jordan blocks of the algebraic monodromy. Let ${\tilde
b}_{f,0}(s)$ be the \emph{microlocal (or reduced)} Bernstein-Sato
polynomial defined by ${b}_{f,0}(s)=(s+1){\tilde b}_{f,0}(s)$.

\medskip

\begin{thm0}\cite{varchenko} \label{varchenko}
Let $M_{f,0}^{n-1}$   be the algebraic
monodromy action on the $(n-1)$-th cohomology $ H^{n-1}(F_{f,0},\bC) $ of the Milnor fibre
of $f$ at the origin.

\begin{enumerate}

\item ${\tilde b}_{f,0}(s)$ is divisible by $(s-\beta)^n$ if and only if $\beta\in (-1,0)$ and $M_{f,0}^{n-1}$ has a Jordan block of size $n$ for the eigenvalue $\exp(2\pi i (\beta))$.

\item  ${\tilde b}_{f,0}(s)$ is divisible by $(s+1-\alpha)^{n-1}$,
with $\alpha\in \bZ_+$,
if and only if $\alpha=0$ and $M_{f,0}^{n-1}$
has a Jordan block of size $n-1$ for the eigenvalue $1$.
\end{enumerate}

\end{thm0}

\medskip




\begin{proof}[\textbf{Proof of Theorem \ref{bernstein}.}] (see introduction)
The proof follows from Theorems \ref{grn} and \ref{varchenko},
together with the fact that for the eigenvalue $\lambda=1$ (resp.
$\lambda \neq 1$), $\Gr_{2n-2}^{W}H^{n-1}(F_{f,0},\bC)_{\lambda}
\neq 0$ if and only if $M_{f,0}^{n-1}$ has a Jordan block of size
$n-1$ (resp. of size $n$), see \cite{vDS} and \cite{st-oslo}.
\end{proof}

\begin{thm0}\label{bernstein-g}
Let $f\in \bc[z_1,\ldots,z_n]$ be a polynomial such that $\f$ has only isolated  singularities. If $s_0=-c(f)$
is a pole of order $n$ of  $\ztop$, then $(s+c(f))^n$ divides the Bernstein-Sato polynomial $b_{f}(s)$.

In such a case there exists an integer $N \geq 1$ such that
$c(f)=1/N$ and either
\begin{itemize}
\item
$N=1$ and $(s+1)^n$ divides $b_{f}(s)$, or
\medskip
\item $N>1$ and
$(s+1/N)^n(s+2/N)^n\dots(s+(N-1)/N)^n(s+1)^n$ divides $b_{f}(s)$.
\end{itemize}
\end{thm0}

\medskip

\begin{proof}
The proof follows from the following two facts. First, a pole of
order $n$ of the local topological zeta function is also a pole of
the global $\ztop$, and conversely, a pole of order $n$ of $\ztop$
is a pole of some $Z^{}_{\operatorname{top},x}(f,s)$ at some point
$x\in \f$. Secondly $b_f(s)$ is the least common multiple of all
local Bernstein-Sato polynomials $b_{f,x}(s)$.
\end{proof}

\medskip

\noindent {\bf 3.2. Non-degenerate Newton polyhedron.} For the
notion of function with non-degenerate Newton polyhedron we refer
for instance to  \cite{lv} or \cite{agv}.  Remark that
almost all polynomials are non-degenerate with respect to their
(either local or global) Newton polyhedron (see \cite{agv} p.151).

For such functions, Denef proved that a set of candidate poles of
the corresponding zeta functions is obtained from the
$(n-1)$-dimensional faces of the corresponding polyhedron, e.g.
see \cite{de:95}. Loeser \cite{lo:90} proved that under some
additional conditions these candidate poles are roots of the
Bernstein-Sato polynomial $b_{f}(s)$.

\medskip

\begin{cor0}
Let $f:(\bc^{n},0)\to(\bc,0)$ be a germ of a holomorphic function
which defines a germ of an isolated hypersurface singularity which
is non-degenerate with respect to its Newton polyhedron at the
origin. If $s_0$ is a pole of order $n$ of  $\ztopl$ then
$(s+s_0)^n$ divides $b_{f,0}(s)$.
\end{cor0}

\medskip

\begin{proof}
Under the hypothesis, Laeremans and the third author proved in
Theorem 2.4 in \cite{lv} that $\ztopl$ has at most one pole of
order $n$. Moreover, if such a pole exists, then it is the pole
closest to the origin which coincides with $-c_0(f)$. Thus after
Theorem \ref{bernstein} we get the result.
\end{proof}


\end{document}